# ON THE EQUALITY OF SYMBOLIC AND ORDINARY POWERS OF BINOMIAL EDGE IDEALS

I. Jahani , Sh. Bayati and F. Rahmati

**ABSTRACT.** In this paper, we investigate whether the symbolic and ordinary powers of a binomial edge ideal $J_G$ are equal. We show that the equality $J_G^t = J_G^{(t)}$ holds for every $t \geq 1$ when $|\operatorname{Ass}(J_G)| = 2$. Moreover, if $G$ is a caterpillar tree, then one has the same equality. Finally, we characterize the generalized caterpillar graphs which the equality of symbolic and ordinary powers of $J_G$ occurs.

## 1 Introduction

Let $G$ be a simple graph on $[n]$, and $S = k[x_1, \ldots, x_n, y_1, \ldots, y_n]$ be the polynomial ring over a field $k$. The binomial edge ideal $J_G$ of $G$ is defined in $S$ as follows:

$$J_G = (x_i y_j - x_j y_i : \{i,j\} \in E(G)),$$

where $E(G)$ is the edge set of $G$. These ideals, first introduced in [8] and [12], have been of great interest in last years and many combinatorial and algebraic properties and invariants of them have been studied; see some of them in the survey [14]. In this paper, we focus on the question that when the equality of symbolic and ordinary powers of binomial edge ideals holds.

Comparison of symbolic and ordinary powers of squarefree monomial ideals corresponding to combinatorial structures are widely studied; see for example [1, 2, 9, 10, 17, 18]. However, less is known about this problem regarding binomial edge ideals. Our aim is to find classes of graphs for which the equality holds for those powers of binomial edge ideals. Ohtani shows in [13, Theorem 4.3] that when $G$ is a complete $r$-partite graph, then $J_G^{(t)} = J_G^t$ for every $t \geq 1$. In Theorem 3.6, we show that for a caterpillar tree $G$, one has $J_G^{(t)} = J_G^t$ for every $t \geq 1$. We also see in Theorem 3.2 that when $|\operatorname{Ass}(J_G)| = 2$, then the equality of symbolic and ordinary powers of $J_G$ holds.

When $G$ is a block graph with Cohen-Macaulay binomial edge ideal, Ene et al. give a complete characterization of those graphs with equality of symbolic and ordinary powers of $J_G$. By dropping the requirement of being Cohen-Macaulay, we give a partial generalization of that result for generalized caterpillar graphs. Although the equality of symbolic and ordinary powers of $J_G$ for a block graph $G$ is equivalent to closedness of $G$ when $J_G$ is Cohen-Macaulay, it turns out in Theorem 3.11 that when $G$ is a generalized caterpillar, that equality holds exactly when $G$ is weakly closed. Moreover, either $G$ is a block graph with Cohen-Macaulay binomial edge ideal or $G$ is a generalized caterpillar graph, then equality of symbolic and ordinary powers is equivalent to $G$ being net-free, see [5, Theorem 4.1] and Theorem 3.11.

Ene and Herzog present a nice result in [4, Theorem 3.3] which shows that if $\text{in}_<(J_G)$, for some monomial order $<$ with a nice property, is normally torsion-free, then $J_G^{(t)} = J_G^t$ for every $t \geq 1$. It turns out that particularly the desired equality holds when $G$ is a closed graph, [4, Corollary 3.4]. One may notice that in this paper, we deal with graphs which are not closed in general.

In Section 2, we see some definition and needed results regarding the symbolic powers of ideals. Next we discuss the comparison of symbolic and ordinary powers in Section 3.

## 2 Preliminaries

In this section, we summarize basic facts about symbolic powers of ideals and binomial edge ideals.

Recall that throughout the paper, $G$ is a simple graph on the vertex set $[n]$ with the edge set $E(G)$, and Suppose that $S = k[x_1, \ldots, x_n, y_1, \ldots, y_n]$ is the polynomial ring over a field $k$ with $2n$ variables. The **binomial edge ideal** $J_G$ of $G$ is the ideal generated by binomials $[i,j] = x_i y_j - y_i x_j$ which $\{i,j\} \in E(G)$.

**Definition 2.1.** *Let $I$ be an ideal of $S$, and suppose that $Min(I)$ is the set of the minimal prime ideals of $I$. For an integer $k \geq 1$, one defines the $k$-th **symbolic power** of $I$ as follows:*

$$I^{(k)} = \bigcap_{p \in Min(I)} (I^k S_P \cap S).$$

Symbolic powers do not, in general, coincide with the ordinary powers. It is known that if $I = Q_1 \cap \ldots \cap Q_m$ is an irredundant primary decomposition of $I$ with $\sqrt{Q_i} = p_i$ for all $i$ and $p_1, \ldots, p_s$ are the minimal prime ideals of $I$, then

$$I^{(k)} = Q_1^{(k)} \cap \ldots \cap Q_s^{(k)}.$$

By [7, Corollary 3.5], for binomial edge ideal $I = J_{G_1}$ and $J = J_{G_2}$, where $G_1$ and $G_2$ are graphs with disjoint sets of vertices, we have the following result.

**Proposition 2.2.** *By the above setting, assume that $I^t \neq I^{t+1}$ and $J^t \neq J^{t+1}$ for $t \leq n - 1$. Then*

$$(I+J)^{(n)} = (I+J)^n$$

*if and only if $I^{(t)} = I^t$ and $J^{(t)} = J^t$ for every $t \leq n$.*

For each subset $U \subseteq [n]$, a prime ideal $P_U(G)$ is defined in [8] as follows: let $T = [n]\setminus U$, and $G_1, \ldots, G_{c(U)}$ be the connected components of the induced subgraph of $G$ on $T$. If for each $i$, $\tilde{G}_i$ is the complete graph on the vertex set $V(G_i)$, then $P_U(G)$ is defined to be

$$P_U(G) = (\{x_i, y_i\}_{i \in U} + J_{\tilde{G}_1} + \ldots + J_{\tilde{G}_{c(U)}}).$$

It is known that $P_U(G)$ is a prime ideal. By [8, Theorem 3.2], $J_G = \bigcap_{U \subset [n]} P_U(G)$.

By [13, Proposition 4.2], the symbolic powers of an ideal of maximal minors of a generic matrix coincide with the ordinary powers. Specially for the complete graph $K_n$ on $[n]$ and by Proposition 2.2, for each $t \geq 1$ we have

$$P_U(G)^{(t)} = P_U(G)^t. \tag{1}$$

See [4] for more details.

In the next proposition by [15, Lemma 4.1], one has a combinatorial equivalence for $|\operatorname{Ass}(J_G)| = 2$.

**Proposition 2.3.** *Let $G$ be a connected graph on $[n]$ and $T_G = \{x \in [n] : deg(x) = n-1\}$. Then $|Ass(J_G)| = 2$ if and only if the following conditions hold:*

1. *$T_G \neq \emptyset$ and $G_{[n]\setminus T_G}$ is disconnected.*
2. *$G_{[n]\setminus T_G}$ is a disjoint union of complete graphs.*

*Furthermore, under the above assumptions, one has $J_G = J_{K_n} \cap P_{T_G}(G)$.*

## 3 Equality of symbolic and ordinary powers

Suppose that $R$ is a polynomial ring over a field $k$ and $\Pi \subseteq R$ is a finite subset with a partial order $\leq$. The ring $R$ is called an **algebra with straightening law** or briefly an ASL (on $\Pi$, over $k$) if the following conditions hold:

1. $\Pi$ consists of homogeneous elements of positive degree and generates $R$ as a $k$-algebra;
2. The products $\xi_1 \ldots \xi_m$, $m \in \mathbb{N}$ and $\xi_i \in \Pi$, such that $\xi_1 \leq \ldots \leq \xi_m$ are linearly independent. They are called **standard monomials**;
3. (Straightening law) For all incomparable $\xi, \nu \in \Pi$ the product $\xi\nu$ has a representation $\xi\nu = \sum a_\mu \mu, a_\mu \in k, a_\mu \neq 0$, with standard monomials $\mu$, satisfying the following condition:
   every $\mu$ contains a factor $\zeta \in \Pi$ such that $\zeta \leq \xi$ and $\zeta \leq \nu$.

**Proposition 3.1.** *([3, Proposition 4.1]) Suppose that $R$ is an ASL over $k$ on $\Pi$. Then:*

1. *The standard monomials generate $R$ as a $k$-module, thus forming a $k$-basis of $R$.*
2. *Furthermore, every monomial $\mu = \xi_1 \ldots \xi_m, \xi_i \in \Pi$ has a standard representation in which every standard monomial contains a factor $\xi \leq \xi_1, \ldots, \xi_m$.*

Suppose that
$$X = \begin{pmatrix} x_1 & \cdots & x_n \\ y_1 & \cdots & y_n \end{pmatrix}$$
is a $2 \times n$ matrix of indeterminants over $k$. Recall that $[i,j]$ is the $2 \times 2$ minor of $X$ by columns $i$ and $j$.

The set of all minors of $X$ can be partially ordered with setting

1. For every $i, j, k$, $x_k < [i,j]$ and $y_k < [i,j]$;
2. If $k \leq k'$, then $x_k \leq x_{k'}, y_{k'}$ and $y_k \leq y_{k'}$;
3. If $i \leq i'$ and $j \leq j'$,
$$[i,j] \leq [i',j'].$$

By setting $\Pi$ to be the set of all minors of $X$, the ring $S$ is an ASL with the above described partial order. For more details, see [3].

**Theorem 3.2.** *Let $G$ be a connected graph on $[n]$ such that $|Ass(J_G)| = 2$. Then for each $t \geq 1$*
$$J_G^t = J_G^{(t)}.$$

*Proof.* We know that $J_G^t \subseteq J_G^{(t)}$. Now we should prove the converse. Since we assume that $|Ass(J_G)| = 2$, by Proposition 2.3, $T_G \neq \emptyset$ and $G_{[n] \setminus T_G}$ is a disjoint union of some complete graphs like $G_1, \ldots, G_c$ and

$$J_G = J_{K_n} \bigcap P_{T_G}(G)$$

where $P_{T_G}(G) = (\{x_i, y_i\}_{i \in T_G}, J_{G_1}, \ldots, J_{G_c})$ and $K_n$ is the complete graph on $[n]$. Since for each $i, j \in \{1, \ldots, c\}$, $V(G_i) \cap V(G_j) = \emptyset$ and $T_G \cap V(G_i) = \emptyset$, so by the Equation (1) and Proposition 2.2, we have

$$J_G^{(t)} = J_{K_n}^t \cap (\{x_i, y_i\}_{i \in T_G}, J_{G_1}, \ldots, J_{G_c})^t.$$

It is sufficient to show that every standard monomial $\mu \in J_G^{(t)}$ is in $J_G^t$. We can consider $\mu$ of the form

$$\mu = N.(\prod_{i=1}^{\alpha}[a_i, b_i])(\prod_{j=1}^{\alpha'}[c_j, d_j])(\prod_{k=1}^{\beta}[e_k, f_k])$$

with the following conditions:

- N is a multiplication of some variables.
- for each $i$, at least one of $a_i$ or $b_i$ is in $T_G$.
- for each $j$, both of $c_j$ and $d_j$ are in $V(G_l)$ for some $l = 1, \ldots, c$.
- for each $k$, there exist distinct numbers $k_0, k_1 \in \{1, \ldots, c\}$ such that $e_k \in V(G_{k_0})$ and $f_k \in V(G_{k_1})$.

In the last case, we consider edges which are not in $E(G)$. Following [13], the minors in the last case are called the bad edges. If $\alpha + \alpha' \geq t$, obviously $\mu \in J_G^t$. So we suppose that $\alpha + \alpha' < t$. We know $\beta > 0$ because $\mu \in J_{K_n}^t$. On the other hand, we know that $\mu \in P_{T_G}(G)^t$. So $\mu$ needs to have at least $t - (\alpha + \alpha')$ variables from the set $\{x_i, y_i\}_{i \in T_G}$. Therefore, $N$ has a factor of the form $x_{i_1}^{\alpha_{i_1}} \ldots x_{i_s}^{\alpha_{i_s}} y_{j_1}^{\alpha'_{j_1}} \ldots y_{j_l}^{\alpha'_{j_l}}$, such that

$$\sum \alpha_{i_l} + \sum \alpha'_{j_r} = t - (\alpha + \alpha').$$

We know that we have at least $t - (\alpha + \alpha')$ bad edges. So for each of them, there exists $x_i$ or $y_j$ with $i, j \in T_G$ and for each bad edge $[a, b]$ of $\mu$, we can obtain it by non-bad edges with the following relation:

$$i_l[a, b] = a[i_l, b] - b[a, i_l]$$

or

$$j_r[a, b] = a[j_r, b] - b[a, j_r].$$

Hence

$$\prod_{k=1}^{\beta}[e_k, f_k] \in J_G^{t-(\alpha+\alpha')}$$

and we conclude that $\mu \in J_G^t$.

□

Suppose that $<_{lex}$ be the lexicographic term order induced by

$$x_1 >_{lex} \ldots >_{lex} x_n >_{lex} y_1 >_{lex} \ldots >_{lex} y_n.$$

Let $i, j \in [n]$ such that $i < j$. Then a path $\pi : i = i_0, i_1, \ldots, i_r = j$ in $G$ from $i$ to $j$ is called **admissible** if the following conditions hold:

1. either $i_k < i$ or $i_k > j$ for every $k = 1, \ldots, r-1$;
2. for each proper subset $\{j_1, \ldots, j_s\}$ of $\{i_1, \ldots, i_{r-1}\}$, the sequence $i, j_1, \ldots, j_s, j$ is not a path.

In particular, all the edges $\{i, j\}$ of $G$, with $i < j$, are admissible paths from $i$ to $j$. Now, let $\pi : i = i_0, i_1, \ldots, i_r = j$ be an admissible path in G. Associated to $\pi$ the following squarefree monomial is defined:

$$u_\pi := \prod_{i_k > j} x_{i_k} \prod_{i_l < i} y_{i_l}.$$

The next theorem provides a reduced Gröbner basis for $J_G$ with respect to the monomial order described above.

**Theorem 3.3.** *([8, Theorem 2.1]) Let $G$ be a graph. Then the following set of binomials in $S$ is a reduced Gröbner basis of $J_G$ with respect to $<_{lex}$ as described above:*

$$\mathcal{G} = \{u_\pi f_{ij} : \pi \text{ is an admissible path from } i \text{ to } j\}.$$

The following result by [4, Theorem 3.3] provides a sufficient condition for a binomial edge ideal to have equal symbolic and ordinary powers by means of its initial ideal.

**Theorem 3.4.** *[4, Theorem 3.3] Let $G$ be a connected graph on the vertex set $[n]$. If $\mathrm{in}_{<_{lex}}(J_G)$ is a normally torsion-free ideal, then $J_G^{(k)} = J_G^k$ for $k \geq 1$.*

Let $\Delta$ be a simplicial complex on the vertex set [n]. A **cycle** or, more precisely, an s-cycle of $\Delta$ ($s \geq 2$) is an alternating sequence of distinct vertices and facets

$$v_1, F_1, \ldots, v_s, F_s, v_{s+1} = v_1$$

such that $v_i, v_{i+1} \in F_i$ for $i = 1, \ldots, s$. The cycle $C$ is called **odd** (**even**) if $s$ is an odd (even) number. A cycle is **special** if it has no facet containing more than two vertices of the cycle.

A simple graph $G$ is called **closed** with respect to a given labelling of the vertices if the following condition is satisfied: For all $\{i, j\}, \{k, l\} \in E(G)$ with $i < j$ and $k < l$, one has $\{j, l\} \in E(G)$ if $i = k$, and $\{i, k\} \in E(G)$ if $j = l$. A simple graph $G$ is closed if there exists a labelling such that it is closed with respect to it.

The notion of an $m$-closed graph is introduced as a generalization of closed graphs in [16]. The graph $G$ is called an $m$-closed graph when its vertices can be labeled by $[n]$ such that the elements of the reduced Gröbner basis of $J_G$ with respect to the above-mentioned lexicographical order have degree at most $m$, and $m$ is the least integer with this property for $G$.

Let $G$ be a graph and $u, v$ be two vertices of $G$. If $v$ is the only vertex of $G$ adjacent to $u$, we say that the edge $e = \{u, v\}$ is a **whisker** of $v$, and $u$ is said to be incident with a whisker of $v$. By adding a whisker to a vertex $v$ of $G$, we mean adding a new vertex $v'$ and an edge $\{v, v'\}$ to $G$.

**Definition 3.5.** *A **caterpillar tree** is a tree $G$ with the property that it contains a path $P$ such that any vertex of $G$ is either a vertex of $P$ or it is adjacent to a vertex of $P$.*

We may assume that the path $P$ in the definition of a caterpillar tree is a longest induced path of $G$ and call it a central path of $G$.

Associated to a monomial ideal $I \subseteq k[x_1, \ldots, x_n, y_1, \ldots, y_n]$ with minimal set of generators $G(I)$, a simplicial complex $\Delta(I)$ is defined whose facets are the sets
$$\{x_{i_1}, \ldots, x_{i_k}, y_{j_1}, \ldots, y_{j_l}\}$$
with $x_{i_1} \ldots x_{i_k} y_{j_1} \ldots y_{i_l} \in G(I)$.

**Theorem 3.6.** *If $G$ is a caterpillar tree and $J_G$ is its binomial edge ideal, then*
$$J_G^t = J_G^{(t)}$$
*for every $t \geq 1$.*

*Proof.* First of all, we fix a central path $P$. Since $P$ is a central path, it has a vertex of degree one, say $v$. We consider a labeling $\sigma : V(G) \to \{1, \ldots, n\}$ which is a one-to-one corresponding with the following properties:

1. $\sigma(v) = 1$.
2. if $w_1, w_2$ are two vertices of $P$ such that $d(v, w_1) < d(v, w_2)$, then
$$\sigma(w_1) < \sigma(z) < \sigma(w_2)$$
   for each vertex $z$ incident with a whisker of $w_1$.

Regarding this labelling on $G$, by [16, Theorem 3.1] we know that the caterpillar trees are 3-closed, so the facets of $\Delta(\text{in}_{<_{lex}}(J_G))$ are of dimension at most 2.

By [9, Corollary 1.6] and [6, Theorem 5.1], in order to show the desired equality of ordinary and symbolic powers of $J_G$, it is sufficient to show that $\Delta(\text{in}_{<_{lex}}(J_G))$ has no special odd cycles.

By the above labeling on $G$, for each $s$ on the central path $P$, the vertices
$$V_{<s} = \{x_i, y_i | i < s\}$$
are not connected to
$$V_{>s} = \{x_i, y_i | i > s\}$$
in $\Delta(\text{in}_{<_{lex}}(J_G))$. More specifically, $\Delta(\text{in}_{<_{lex}}(J_G)) - \{y_s\}$ is not connected and every cycle in $G$ that intersects both $V_{<s}$ and $V_{>s}$ must include $y_s$ twice and hence it is not special. As a result, we can only have special cycles on induced subcomplex on $V_{<s}$ or $V_{>s}$ for each vertex $s$ on the central path $P$. So we can reduce the problem to $G = K_{1,t}$.

We will prove this by induction on $t$. It is obvious for $t = 1$. Now suppose that for $t = k - 1$, the simplicial complex $\Delta(\text{in}_{<_{lex}}(J_{K_{1,k-1}}))$ has no special odd cycles. Let $t = k$ and $H = \Delta(\text{in}_{<_{lex}}(J_{K_{1,k}}))$.

One can find the facets of $H$ by Theorem 3.3. Suppose that $G_j$ be facets of the form $\{x_1, y_j\}$ for $2 \leq j \leq k$ and $F_{ij}$ be facets of the form $\{y_1, x_i, y_j\}$ for $2 \leq i < j \leq k$. Notice that the facets $G_j$'s and $F_{ij}$'s are the only facets of $H$. So $x_k$ divides none of the minimal monomial generators of the ideal $\text{in}_{<_{lex}}(J_{K_{1,k}})$. So if $H$ has a special odd cycle $C$, it should include $y_k$; Otherwise a contradiction to the induction hypothesis.

We first claim that there exists no special odd cycle in $H$ containing $y_1$. Assume that $y_1$ appears in $C$. Regarding the facets and vertices appearing after $y_1$ in a cycle $C$, we have one of the following cases:

1. First let $y_1, F_{ij}, x_i$ appears in $C$ where $2 \leq i$. Since each facet including $x_i$ with $i \geq 2$ has $y_1$ as well, we can not proceed to complete $y_1, F_{ij}, x_i$ to a special cycle.

2. Next consider the subsequence $y_1, F_{ij}, y_j$ appears in $C$. Since each facet $F_{rj}$ including $y_j$ has $y_1$, we can only proceed to $y_1, F_{ij}, y_j, G_j, x_1$ in $C$. Furthermore, $x_1$ is only a vertex of the facets $G_s$, with $2 \leq s \leq k$. If $s \neq k$, the subsequence $y_1, F_{ij}, y_j, G_j, x_1, G_s, y_s$ can not complete to form a special sequence, because except than $G_s$, each facet including $y_s$ has $y_1$ as well.

Hence regarding the subsequence that appears after $y_1$ in $C$ we can only have either

$$y_1, F_{ij}, y_j, G_j, x_1, G_k, y_k$$

or

$$y_1, F_{ik}, y_k.$$

By symmetry, regarding the subsequence appearing before $y_1$ in $C$, either the subsequence

$$y_k, G_k, x_1, G_{j'}, y_{j'}, F_{ij'}, y_1$$

or

$$y_k, F_{i'k}, y_1$$

in $C$. So in any case, we have an odd number of facets after $y_1$ in $C$, and an odd number of facets before that. So there is no odd cycle including $y_1$.

Hence to construct a special odd cycle, we can only have vertices $x_i$, $y_j$ of $F_{ij}$'s and $x_1$, $y_j$ of $G_j$'s. In particular, notice that from each facet, we have one vertex of type $x_i$'s and one vertex of type $y_j$'s in $C$. So $x_i$'s exactly appear every other steps, as well as $y_j$'s. Thus, the length of the cycle $C$ could be even, a contradiction.

□

Matsuda in [11] introduced the notion of a weakly closed graph as a generalization of closed graphs.

**Definition 3.7.** *The graph $G$ is said to be **weakly closed** if there exists a labelling of the vertices which satisfies the following condition: for all integers $1 \leq i < j < k \leq n$, if $\{i, k\} \in E(G)$, then $\{i, j\} \in E(G)$ or $\{j, k\} \in E(G)$.*

**Definition 3.8.** *Let $P = ([n], <_P)$ be a partially ordered set. The graph $G(P)$ associated to $P$ is a graph on the vertex set $[n]$ such that $\{i, j\} \in E(G(P))$ with $i < j$ if $i <_P j$. A graph $G$ is **comparability** if there exists a partially ordered set $P$ such that $G = G(P)$.*

We have the next theorem from [11, Theorem 1.9] that illustrates the relation between weakly closedness of a graph $G$ and comparability of its complement $\bar{G}$.

**Theorem 3.9.** *Let $G$ be a graph. Then the following assertions are equivalent:*

1. *$G$ is weakly closed.*
2. *$G$ is co-comparability, i.e. its complement graph $\bar{G}$ is comparability.*

A cutpoint of a simple graph $G$ is a vertex $v \in V(G)$ for which the subgraph $G - v$ is disconnected. A clique of the graph $G$ is a complete subgraph of $G$. A block of $G$ is a connected subgraph of $G$ that has no cutpoints and is maximal with respect to this property. A graph $G$ is called a block graph if each block of G is a clique.

Now we introduce a generalization of caterpillar trees:

Let $G$ be a graph on $V$ and $e = \{v_1, v_2\}$ be an edge of $G$. By a **clique join** on $G$ via $e$ we mean the graph obtained by attaching a complete graph $K_t$ to $G$ on $e$ for some $t \geq 2$. More precisely, this graph denoted by $G \sqcup_e K_t$ is the graph whose set of vertices is obtained by adding $t-2$ new vertices $u_1, \ldots, u_{t-2}$ to $V$, and $f$ is an edge of $G \sqcup_e K_t$ if $f$ is either an edge of $G$ or $f = \{u_i, v_j\}$ for some $i = 1, \ldots, i_{t-2}$ and $j = 1, 2$.

Let $G$ be a graph. We call $G$ a **generalized caterpillar graph** if there exists a caterpillar tree $H$ with a central path $P$ of $H$ such that $G$ can be obtained by clique join of some complete graphs $K_{t_1}, \ldots, K_{t_m}$ in succession on $H$ via pairwise distinct edges $e_1, \ldots, e_m$ of $P$ and finally, possibly adding an arbitrary number of whiskers to some vertices. In this definition, we call $P$ a central path of $G$ if among different possible choices of caterpillar tree $H$, the path $P$ is a longest one.

Figure 1 is an example of a generalized caterpillar graph. Notice that each generalized caterpillar graph is a block graph.

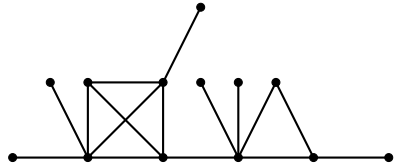

Figure 1

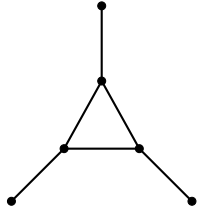

Figure 2

A graph is called **net-free** if it does not have an induced subgraph as Figure 2.

**Lemma 3.10.** *A generalized caterpillar graph is weakly closed if and only if it is net-free.*

*Proof.* First, suppose that $G$ is net-free. We give a labeling on $G$ which walking on a fixed central path of $G$, in each step we first give a label to the vertex of the path, next to vertices incident with its whiskers, and finally to adjacent cliques and their whiskers' vertices. More specifically, according to the definition of generalized caterpillar graphs we may assume that for every two adjacent vertices $v_i$ and $v_j$ in a central path there exists a clique that contains both of them as vertices. We denote this clique by $K_{i,j}$. In central path $P$, there exists a vertex included in exactly one maximal clique, say $v$. We set a labeling on a generalized caterpillar graph by a one-to-one corresponding $\sigma : V(G) \to \{1, \ldots, n\}$ which has the following properties:

1. $\sigma(v) = 1$
2. if $v_i, v_j$ are two vertices of $P$ such that $d(v, v_i) < d(v, v_j)$, then
$$\sigma(v_i) < \sigma(v_j)$$

3. Let $v_i$ and $v_j$ be vertices of $P$ with $\sigma(v_i) < \sigma(v_j)$ and suppose that $w$ is adjacent to $v_i$ by a whisker. Assume that either $z$ is a vertex of $K_{ij}$ other than $v_i$ and $v_j$ or $z$ is adjacent to $K_{ij}$ by a whisker but not adjacent to $P$. Then we set

$$\sigma(v_i) < \sigma(w) < \sigma(z) < \sigma(v_j).$$

With the above labeling on vertices of $G$, we can conclude that $G$ is weakly closed.

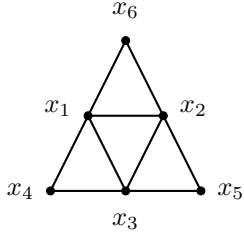

Figure 3

Conversely, suppose that $G$ contains a net as an induced subgraph. Then the complement $\bar{G}$ of $G$ contains a graph in Figure 1 as an induced subgraph. By Theorem 3.9, it is enough to show that $\bar{G}$ (Figure 3) is not comparability. Assume that $\bar{G}$ is comparability. Without loss of generality, let $x_1 < x_2 < x_3$.

Now if $x_2 > x_5$ then $x_2 > x_6$ because $\{x_5, x_6\}$ is not an edge of $\bar{G}$. So

$$x_3 > x_2 > x_6,$$

which is a contradiction because $\{x_3, x_6\}$ is not an edge of $\bar{G}$. Hence $x_1 < x_2 < x_5$, but this is again, a contradiction to the fact that $\{x_1, x_5\} \notin E(\bar{G})$. $\square$

For block graphs with Cohen-Macaulay binomial edge ideals, those one with equality of ordinary and symbolic powers are characterized in [5, Theorem 4.1]. We give a partial generalization of that result.

**Theorem 3.11.** *Suppose that $G$ is a generalized caterpillar graph. Then the following conditions are equivalent:*

1. *$G$ is weakly closed;*
2. *$J_G^i = J_G^{(i)}$ for every $i \geq 1$;*
3. *$J_G^i = J_G^{(i)}$ for some $i \geq 2$;*
4. *$J_G^2 = J_G^{(2)}$;*
5. *$G$ is net-free.*

*Proof.* By Lemma 3.10 we know that the conditions (1) and (5) are equivalent. The conclusions (2) $\Rightarrow$ (3), (2) $\Rightarrow$ (4) and (4) $\Rightarrow$ (3) are trivial. The implications (3) $\Rightarrow$ (5) and (4) $\Rightarrow$ (5) are proved in general case in [5, Theorem 4.1].

(5) $\Rightarrow$ (2) Now we show that if $G$ is a net-free graph, then the symbolic and ordinary powers of $J_G$ are equal. If $G$ is a net-free graph, each 3-vertex clique $D$ of $G$ can have whiskers on at most 2 vertices. Let assume that $D$ has 2 vertices with whiskers. The third vertex of $D$ can not

be the first (or last) vertex of $P$, because we consider the longest path as the central path, so one of the whiskers will be the first (or last) edge of $P$. On the other hand, since $G$ is a net-free graph, these 2 vertices of $D$ which have whiskers should appear on the central path $P$.

By the labeling of $G$ described in Lemma 3.10, for each $s$ on the central path $P$, the vertices
$$V_{<s} = \{x_i, y_i | i < s\}$$
in $\Delta(\text{in}_{<_{lex}}(J_G))$ are not connected to
$$V_{>s} = \{x_i, y_i | i > s\}.$$
More specifically, $\Delta(\text{in}_{<_{lex}}(J_G)) - \{y_s\}$ is not connected and every cycle in $G$ that intersects both $V_{<s}$ and $V_{>s}$ must include $y_s$ twice and hence it is not special. As a result, we can only have special cycles on induced subcomplex on $V_{<s}$ or $V_{>s}$ for each vertex $s$ on the central path $P$, so we may assume that $G$ is a complete graph with some whiskers. On the other hand, notice that by the above discussion if a 3-vertex clique has 2 vertices with whiskers, then the vertices with whiskers lie on the central path $P$. So we can reduce the problem to the case that $G$ is a complete graph with whiskers on at most one vertex. Hence we may assume that $G$ is a complete graph $K_{t'}$ with $t$ whiskers on at most one of its vertices.

If $t = 0$, we have a complete graph. Each complete graph is a closed graph, so in this case we have the desired result by [4, Corollary 3.4]. Hence we can assume that $t \neq 0$. We will prove the rest by induction on $t'$. For $t' = 1, 2$, the graph $G$ is a caterpillar tree. So by Theorem 3.6, we know that $\Delta(\text{in}_{<_{lex}}(J_G))$ has no special odd cycles. Let $t' = k > 2$ and $H = \Delta(\text{in}_{<_{lex}}(J_G))$.

According to the labeling of $G$ described in Lemma 3.10, the vertices of $K_k$ have the labels $1, t+2, t+3, \ldots, t+k$, and the vertex 1 has the whiskers which incident with vertices labeled by $2, 3, \ldots, t+1$. By Theorem 3.3 we can find the facets of $H$. Suppose that $G_{ij}$ is the facet $\{x_i, y_j\}$ for $i = 1$ and $2 \leq j \leq t+k$ or $t+2 \leq i < j \leq t+k$ and $F_{ij}$ is the facets $\{y_1, x_i, y_j\}$ for $2 \leq i \leq t+1$, $i < j \leq t+k$. Notice that the facets $G_{ij}$'s and $F_{ij}$'s are the only facets of $H$. So $x_{t+k}$ divides none of the minimal monomial generators of the ideal $\text{in}_{<_{lex}}(J_G)$. If $H$ has a special odd cycle $C$, it should include $y_{t+k}$; Otherwise a contradiction to the induction hypothesis.

We claim that if $C$ is a special cycle of $H$ including $y_1$, then it is a special even cycle.
Assume that $y_1$ appears in $C$. Regarding the facets and vertices appearing after $y_1$ in $C$, we have one of the following cases:

1. First let the subsequence $y_1, F_{ij}, x_i$ appear in $C$ where $2 \leq i \leq t+1$. Since $x_i$ with $2 \leq i \leq t+1$ does not appear in $G_{ij}$'s and each facet including such $x_i$'s has $y_1$ as well, we can not proceed to complete $y_1, F_{ij}, x_i$ to a special cycle.

2. Next consider the subsequence $y_1, F_{ij}, y_j$ appears in $C$. Since each facet $F_{rj}$ including $y_j$ has $y_1$, we can only proceed to $y_1, F_{ij}, y_j, G_{sj}, x_s$ in $C$. Furthermore, $x_s$ is only a vertex of the facets $G_{s\ell}$'s, with $2 \leq \ell \leq t+k$. Since except than $G_{s\ell}$'s, each facet including $y_\ell$ has $y_1$ as well, the subsequence $y_1, F_{ij}, y_j, G_{sj}, x_s, G_{s\ell}, y_\ell$ can continue only with $G_{pq}$'s to form a special sequence.

   Hence regarding the subsequence that appears after $y_1, F_{ij}, y_j$ in $C$ we can only have
   $$y_1, F_{ij}, y_j, G_{sj}, x_s, G_{sh}, \ldots, G_{rk}, y_{t+k},$$
   with probably no $G_{pq}$'s.

So the only special cycle which includes $y_1$ is of the form of

$$y_{t+k}, G_{u\,t+k}, \ldots, x_{s'}, G_{s'j'}, y_{j'}, F_{ij'}, y_1, F_{ij}, y_j, G_{sj}, x_s, \ldots, G_{v\,t+k}, y_{t+k}$$

which is an even cycle.

Hence to construct a special cycle, we can only have vertices $x_i$, $y_j$ of $F_{ij}$'s and $x_i$, $y_j$ of $G_{ij}$'s. In particular, notice that from each facet, we have one vertex of type $x_i$'s and one vertex of type $y_j$'s in $C$. So $x_i$'s and $y_j$'s exactly appear every other steps, as well as $y_j$'s. Thus, the length of the cycle $C$ could be even, a contradiction.

□

# References

[1] S. Bayati, F. Rahmati. *Squarefree vertex cover algebras.* Communications in Algebra, 42(4), pp. 1518–1538, 2014.

[2] C. Bahiano. *Symbolic powers of edge ideals.* Journal of Algebra, 273, pp. 517–537, 2004.

[3] W. Bruns, U. Vetter. *Determinantal rings*, volume 1327. Springer, 2006.

[4] V. Ene, J. Herzog. *On the symbolic powers of binomial edge ideals.* Combinatorial Structures in Algebra and Geometry, pp. 43–50, 2018.

[5] V. Ene, G. Rinaldo, N. Terai. *Powers of binomial edge ideals with quadratic gröbner bases.* Nagoya Mathematical Journal, pp. 1–23, 2021.

[6] D. Fulkerson, A. Hoffman, R. Oppenheim. *On balanced matrices.* Pivoting and Extension, pp. 120–132. Springer, 1974.

[7] H. T. Hà, H. D. Nguyen, N. V. Trung, T. N. Trung. *Symbolic powers of sums of ideals.* Mathematische Zeitschrift, 294(3), pp. 1499–1520, 2020.

[8] J. Herzog, T. Hibi, F. Hreinsdóttir, T. Kahle, J Rauh. *Binomial edge ideals and conditional independence statements.* Advances in Applied Mathematics, 45(3), pp. 317–333, 2010.

[9] J. Herzog, T. Hibi, N. Trung, X. Zheng. *Standard graded vertex cover algebras, cycles and leaves.* Transactions of the American Mathematical Society, 360(12), pp. 6231–6249, 2008.

[10] J. Herzog, T. Hibi, N. V. Trung. *Symbolic powers of monomial ideals and vertex cover algebras.* Advances in Mathematics, 210(1), pp. 304–322, 2007.

[11] K. Matsuda. *Weakly closed graphs and f-purity of binomial edge ideals.* Algebra Colloquium, 25, pp. 567–578, 2018.

[12] M. Ohtani. *Graphs and ideals generated by some 2-minors.* Communications in Algebra, 39(3), pp. 905–917, 2011.

[13] M. Ohtani. *Binomial edge ideals of complete multipartite graphs.* Communications in Algebra, 41(10), pp. 3858–3867, 2013.

[14] S. Saeedi Madani. *Binomial edge ideals: A survey.* The 24th National School on Algebra, pp. 83–94. Springer, 2016.

[15] L. Sharifan. *Binomial edge ideals with special set of associated primes.* Communications in Algebra, 43(2), pp. 503–520, 2015.

[16] L. Sharifan, M. Javanbakht. *On m-closed graphs.* arXiv preprint, arXiv:1708.08864, 2017.

[17] S. Sullivant. *Combinatorial symbolic powers.* Journal of Algebra, 319(1), pp. 115–142, 2008.

[18] N. V. Trung, T. M. Tuan. *Equality of ordinary and symbolic powers of stanley-reisner ideals.* Journal of Algebra, 328(1), pp. 77–93, 2011.

Department of Mathematics and Computer Science, Amirkabir University of Technology (Tehran Polytechnic), Iran

*E-mail address*: `imanjahani@aut.ac.ir`
*E-mail address*: `shamilabayati@gmail.com`
*E-mail address*: `frahmati@aut.ac.ir`